\documentclass[article,11pt]{amsart}

\usepackage{amsthm}
\usepackage{amssymb}
\usepackage{amstext}
\usepackage{amsmath}
\usepackage{url}
\usepackage[colorlinks]{hyperref}
\hypersetup{
linkcolor=blue,
citecolor=blue,
}
\usepackage[margin = 40pt,font=small,labelfont=bf]{caption}
\usepackage{arydshln}
\usepackage{pdfsync}
\usepackage[mac]{inputenc}
\usepackage[T1]{fontenc}
\usepackage{color}
\usepackage{algorithm}
\usepackage{algorithmic}
\usepackage{enumerate}
\usepackage{bbm}
\usepackage{ae,aecompl}
\usepackage{pdfpages}
\usepackage{epstopdf}
\usepackage{graphicx}
\usepackage{epsfig}
\usepackage{subfigure}

\numberwithin{equation}{section}
\theoremstyle{plain}

\def\build#1_#2^#3{\mathrel{\mathop{\kern 0pt#1}\limits_{#2}^{#3}}}

\def\videbox{\mathbin{\vbox{\hrule\hbox{\vrule height1.4ex \kern.6em\vrule height1.4ex}\hrule}}}

\newcommand{\bv}{\boldsymbol{v}}

\newcommand{\rI}{\mathrm{I}}

\newcommand{\E}{{\mathbb E}}
\newcommand{\R}{{\mathbb R}}

\newcommand{\HH}{\ensuremath{\mathcal H}}
\newcommand{\XX}{\ensuremath{\mathcal X}}

\def\diag{\mathop{\rm diag}\nolimits}%

\newcommand{\thefont}[2]{\fontsize{#1}{#2}\fontshape{n}\selectfont}
\newcommand{\1}{\rlap{\thefont{10pt}{12pt}1}\kern.16em\rlap{\thefont{11pt}{13.2pt}1}\kern.4em}

\begin{document}
\title[Stochastic algorithms for optimal transport]
{Erratum : ``Asymptotic distribution  and convergence rates of stochastic algorithms for entropic optimal transportation between probability measures" \vspace{1ex}}
\author{Bernard Bercu and J\'{e}r\'{e}mie Bigot}
\dedicatory{\normalsize Universit\'e de Bordeaux\\
Institut de Math\'ematiques de Bordeaux et CNRS  (UMR 5251)}

\maketitle


This erratum offers a correction of paper ``Asymptotic distribution  and convergence rates of stochastic algorithms for entropic optimal transportation between probability measures" by B. Bercu and J. Bigot posted on Arxiv at: \url{https://arxiv.org/abs/1812.09150v5}, and that has been published in {\it Annals of Statistics 49(2): 968-987 (2021)} \cite{BB21}. \\

In Appendix A of  paper \cite{BB21}, inequality (A.4) in Lemma A.1 states that
\begin{equation}
\rho_{A_\varepsilon}(v^{\ast})  \geq \frac{1}{\varepsilon} \min_{1 \leq j \leq J} \nu_{j}.
\end{equation}
However, this lower bound for   $\rho_{A_\varepsilon}(v^{\ast})$ is not correct. Using arguments in the proof of Lemma A.1, the right expression for this lower bound is

\setcounter{equation}{3}
\begin{equation} \label{eq:rightbound}
\rho_{A_\varepsilon}(v^{\ast})  \geq  \frac{1}{\varepsilon} \E \Big[   \min_{1 \leq j \leq J}   \pi_{j}(X,v^{\ast})  \Big]
\end{equation}
where $v^{\ast}$ is the optimal dual Kantorovich potential and $\pi_j(X,v^{\ast})$ defined by (3.3) in paper  \cite{BB21}, satisfies for all $1 \leq j \leq J$, $\nu_j =\E \bigl[ \pi_j(X,v^{\ast}) \bigr]$.
\vspace{1ex} \\
We thank Lenaïc Chizat and Alex Delalande for pointing this out to us. Changing this lower bound does not affect the results of the paper as they have been obtained for a fixed value of $\varepsilon >0$. Moreover, the proofs of the paper  do not use the lower bound (A.4). 
\vspace{1ex} \\
\begin{center}
\bf{PROOF OF THE LOWER BOUND IN LEMMA A.1}
\end{center}
\ \\
\setcounter{section}{1}
\setcounter{equation}{0}
For any $v \in \R^J$, let $A_\varepsilon(v)$ be the negative semi-definite matrix defined by
$$
 A_\varepsilon(v)= -\frac{1}{\varepsilon} \E \bigl[ A_{\varepsilon}(X,v)  \bigr] 
$$
where, for all $x \in \XX$,
$$
A_{\varepsilon}(x,v) =  \diag(\pi(x,v))  - \pi(x,v) \pi(x,v)^T,
$$
with $\pi_j(x,v)$ defined by (3.3) in paper  \cite{BB21}.
We denote by  $\lambda_{1}^\varepsilon(v), \ldots, \lambda_{J}^\varepsilon(v)$  the real eigenvalues of the matrix $A_\varepsilon(v)$ that is of rank $J-1$ with $\lambda_{J}^\varepsilon(v) = 0$ 
and $\lambda_j^\varepsilon(v) < 0$ for all $1 \leq j \leq J-1$. We also recall that $\bv_{J}  =\frac{1}{\sqrt{J}}\mathbf{1}_J$ is the unit eigenvector associated to $\lambda_{J}^\varepsilon(v)$
and we denote by $\langle \bv_J \rangle$ the one-dimensional subspace of $ \R^J$  spanned by $\bv_J$. Hereafter, choose $v=v^\ast$ and let
$$
\rho_{A_\varepsilon}(v^{\ast})=-\max_{1\leq j \leq J-1} \big\{ \lambda_{j}^\varepsilon(v^\ast) \big\}=\min_{1\leq j \leq J-1} \big\{ -\lambda_{j}^\varepsilon(v^\ast) \big\}.
$$
We clearly have
$$
\rho_{A_\varepsilon}(v^{\ast})= \min_{ \left\{ u \in \langle \bv_J \rangle^\perp \;  : \; \| u \| = 1\right\}} -u^T A_\varepsilon(v^{\ast}) u = \frac{1}{\varepsilon} \min_{ \left\{ u \in \langle \bv_J \rangle^\perp \;  : \; \| u \| = 1\right\}} u^T \E \bigl[ A_{\varepsilon}(X,v^\ast)  \bigr] u.
$$
However, it follows from Theorem 6 and inequality (5.11) in \cite{Steerneman:2005:PM} that the second smallest eigenvalue of $A_{\varepsilon}(x,v^{\ast})$ 
is lower bounded by
$
\min_{1 \leq j \leq J} \pi_{j}(x,v^{\ast})
$
for all $x \in \XX$. Moreover, $A_{\varepsilon}(x,v^{\ast})$  is a positive semi-definite matrix  of rank $J-1$, and its eigenvector associated to its smallest eigenvalue (equal to  zero)  is  $\bv_J$ for any $x \in \XX$. Consequently, for all unit vector $u \in \langle \bv_J \rangle^\perp$,
$$
u^T \E \big[ A_{\varepsilon}(X,v^\ast)  \big] u \geq \E \big[ \min_{1 \leq j \leq J} \pi_{j}(X,v^{\ast})\big]
$$
which leads to 
\begin{equation*}
\rho_{A_\varepsilon}(v^{\ast})  \geq  \frac{1}{\varepsilon} \E \big[   \min_{1 \leq j \leq J}   \pi_{j}(X,v^{\ast})  \big].
\end{equation*}
\vspace{1ex} \\
\begin{center}
\bf{ILLUSTRATIVE EXAMPLE}
\end{center}
\ \\
We conclude this erratum by providing an illustrative example suggested by Alex Delalande for which one may obtain  explicit expressions of both $\rho_{A_\varepsilon}(v^{\ast})$ and the lower bound \eqref{eq:rightbound} as functions of $\varepsilon$ that allows to study their convergence as  $\varepsilon$ goes to zero. 
Assume that the dimension $d=1$, that $\mu$ is the uniform distribution on $[-1/2,1/2]$ and that
$$
\nu=\frac{1}{2} \delta_{y_1}+\frac{1}{2} \delta_{y_2} \quad \mbox{with} \quad y_1=-1, y_2=1.
$$
Since $\mu(]-\infty,0]) = \mu([0,+\infty[)  =  \frac{1}{2}$, it follows  by symmetry of $\mu$ and $\nu$ that  
$$
v^\ast= \begin{pmatrix}
0 \\
0
\end{pmatrix}.
$$
Hence, if one considers  the cost function $c(x,y) = -xy$, one obtains that
$$
\pi_1(x,v^{\ast}) = \frac{e^{-x/\varepsilon}}{e^{-x/\varepsilon}+e^{x/\varepsilon}} \hspace{1cm} \mbox{and} \hspace{1cm} \pi_2(x,v^{\ast})  = \frac{e^{x/\varepsilon}}{e^{-x/\varepsilon}+e^{x/\varepsilon}}.
$$
Consequently,
$$
A_{\varepsilon}(x,v^{\ast})=\begin{pmatrix}
\pi_1(x,v^{\ast})(1-\pi_1(x,v^{\ast})) & -\pi_1(x,v^{\ast})\pi_2(x,v^{\ast})  \\
-\pi_1(x,v^{\ast})\pi_2(x,v^{\ast}) & \pi_2(x,v^{\ast})(1-\pi_2(x,v^{\ast}))
\end{pmatrix}
$$
Hence, we obtain that
$$
\E \big[ A_{\varepsilon}(X,v^\ast)  \big]= \varepsilon a_\varepsilon
\begin{pmatrix}
1 &-1 \\
-1 & 1
\end{pmatrix}
$$
where
\begin{align*}
a_\varepsilon &= \frac{1}{\varepsilon}\int_{-1/2}^{1/2} \frac{1}{(e^{-x/\varepsilon}+e^{x/\varepsilon})^2}\,dx=
\int_{-1/2\varepsilon}^{1/2\varepsilon} \frac{1}{(e^{-y}+e^{y})^2}\,dy \\
&= -\frac{1}{2} \left[ \frac{1}{1+e^{2y}} \right]_{-1/2\varepsilon}^{1/2\varepsilon}=
-\frac{1}{2}\left( \frac{1}{1+e^{1/\varepsilon}} -  \frac{1}{1+e^{-1/\varepsilon}} \right).
\end{align*}
Therefore, we find that
$$
A_\varepsilon(v^\ast) = - a_\varepsilon
\begin{pmatrix}
1 &-1 \\
-1 & 1
\end{pmatrix}
$$
and $\rho_{A_\varepsilon}(v^{\ast})=2 a_\varepsilon$. In addition,
$$
\lim_{\varepsilon \rightarrow 0} a_\varepsilon=\frac{1}{2} \hspace{1cm} \mbox{and} \hspace{1cm} 
\lim_{\varepsilon \rightarrow 0} \rho_{A_\varepsilon}(v^{\ast})=1.
$$
Finally, we have
$$
\E \big[\min(\pi_{1}(X,v^{\ast}), \pi_{2}(X,v^{\ast})) \big]=\E \big[ \pi_{2}(X,v^{\ast}) \rI_{X<0} + \pi_{1}(X,v^{\ast}) \rI_{X>0}\big]=\varepsilon b_\varepsilon
$$
where
\begin{align*}
b_\varepsilon &= \frac{1}{\varepsilon}\left(\int_{-1/2}^{0} \frac{e^{x/\varepsilon}}{e^{-x/\varepsilon}+e^{x/\varepsilon}}\,dx+\int_{0}^{1/2} \frac{e^{-x/\varepsilon}}{e^{-x/\varepsilon}+e^{x/\varepsilon}}\,dx \right) 
\\
&= \int_{-1/2\varepsilon}^0 \frac{e^{y}}{e^{-y}+e^{y}}\,dy +\int_{0}^{1/2\varepsilon} \frac{e^{-y}}{e^{-y}+e^{y}}\,dy \\
&=  \Big[ \frac{1}{2}\log(1+e^{2y}) \Big]_{-1/2\varepsilon}^0+ \Big[ y - \frac{1}{2}\log(1+e^{2y}) \Big]_{0}^{1/2\varepsilon}\\
&=\log 2-\log(1+e^{-1/\varepsilon}).
\end{align*}
We have
$$
\lim_{\varepsilon \rightarrow 0} b_\varepsilon=\log 2.
$$
Moreover,
$$
\rho_{A_\varepsilon}(v^{\ast}) \geq \E \big[\min(\pi_{1}(X,v^{\ast}), \pi_{2}(X,v^{\ast})) \big] \Longleftrightarrow 2 a_\varepsilon \geq b_\varepsilon
$$
which is of course true for all $\varepsilon >0$.
Therefore, the above calculation show that $\rho_{A_\varepsilon}(v^{\ast})$ and the lower bound  \eqref{eq:rightbound} remain finite as  $\varepsilon$ goes to zero. These findings are in agreement with recent research works \cite{pmlr-v151-delalande22a}[Theorem 3.2] and \cite{chizat2024}[Proposition 5.1] on the computation of lower bounds for $\rho_{A_\varepsilon}(v^{\ast})$ in the semi-discrete setting when $\mu$ is an absolutely continuous probability measure with an upper and lower bounded density $f_{\mu}$. 
Using the fact that $\lim_{\varepsilon \to 0} \rho_{A_\varepsilon}(v^{\ast}) = 1$, we can guess from the above computations that the un-regularized semi-dual OT functional $H_0(v)$, defined in Equation (2.10) of  paper \cite{BB21}, should be twice differentiable at its optimiser $v^\ast$ with Hessian given by
$$
A_0(v^\ast) = -\frac{1}{2}
\begin{pmatrix}
1 & -1\\
-1 & 1
\end{pmatrix}.
$$
To confirm this guess, we can use \cite{newton19}[Theorem 1.3] with $d=1$ that ensures that, for $i \neq j \in \{1,2\}$,
$$
\frac{\partial^2}{\partial v_i \partial v_j} H_0(v^\ast) = \int_{T^{-1}(y_i) \cap T^{-1}(y_j) } \frac{f_{\mu}(x)}{|\frac{\partial}{\partial x} c(x,y_i) - \frac{\partial}{\partial x} c(x,y_j)|} d \HH^{d-1}(x),
$$
where $T$ is the unique optimal mapping from $\mu$ to $\nu$, 
$T^{-1}(y_i)$ and $T^{-1}(y_j)$ are the so-called Laguerre cells, and $\HH^{d-1}$ denotes the $(d-1)$-dimensional Hausdorff measure. For this illustrative example, one has that
$$
T^{-1}(-1) = ]-\infty,0] \quad \mbox{and} \quad T^{-1}(1) = [0,+\infty[.
$$
Therefore, one has that
$T^{-1}(y_i) \cap T^{-1}(y_j) = \{ 0 \}$. Moreover, as $\frac{\partial}{\partial x} c(x,y) = -y$, we obtain that, for $i \neq j \in \{1,2\}$,
$$
\left|\frac{\partial}{\partial x} c(x,y_i) - \frac{\partial}{\partial x} c(x,y_j)\right| = |1-(-1)| = 2.
$$
Hence,  since $f_{\mu}(0) =1$, $\HH^0(\{0\}) = 1$, we have that, for $i \neq j$,
$$
\frac{\partial^2}{\partial v_i \partial v_j} H_0(v^\ast) = \frac{1}{2}
$$
from which it follows that
$$
\nabla^2 H_0(v^\ast) = A_0(v^\ast) = -\frac{1}{2}
\begin{pmatrix}
1 & -1\\
-1 & 1
\end{pmatrix}.$$

\bibliographystyle{acm}
\bibliography{AlgoStoOT-Erratum-Long}

\begin{thebibliography}{1}

\bibitem{BB21}
{\sc Bercu, B., and Bigot, J.}
\newblock {Asymptotic distribution and convergence rates of stochastic
  algorithms for entropic optimal transportation between probability measures}.
\newblock {\em The Annals of Statistics 49}, 2 (2021), 968 -- 987.

\bibitem{chizat2024}
{\sc Chizat, L., Delalande, A., and Vaskevicttius, T.}
\newblock Sharper exponential convergence rates for sinkhorn's algorithm in
  continuous settings, 2024.

\bibitem{pmlr-v151-delalande22a}
{\sc Delalande, A.}
\newblock Nearly tight convergence bounds for semi-discrete entropic optimal
  transport.
\newblock In {\em Proceedings of The 25th International Conference on
  Artificial Intelligence and Statistics\/} (28--30 Mar 2022), G.~Camps-Valls,
  F.~J.~R. Ruiz, and I.~Valera, Eds., vol.~151 of {\em Proceedings of Machine
  Learning Research}, PMLR, pp.~1619--1642.

\bibitem{newton19}
{\sc Kitagawa, J., Mérigot, Q., and Thibert, B.}
\newblock Convergence of a newton algorithm for semi-discrete optimal
  transport.
\newblock {\em Journal of the European Mathematical Society 21}, 9 (2019).

\bibitem{Steerneman:2005:PM}
{\sc Steerneman, T., and {van Perlo-ten Kleij}, F.}
\newblock Properties of the matrix {$A - X Y*$}.
\newblock {\em Linear Algebra and its Applications 410}, 1 (2005), 70--86.

\end{thebibliography}

\end{document}